\documentclass[11pt]{amsart}

\usepackage{amssymb, amscd, latexsym,color,amscd}
\usepackage{url}
\usepackage[final]{graphicx}

\newtheorem{thm}{Theorem}[section]
\newtheorem{prop}[thm]{Proposition}
\newtheorem{lem}[thm]{Lemma}
\newtheorem{cor}[thm]{Corollary}

\newtheorem{problem}[thm]{Problem}

\newtheorem{construction}[thm]{Construction}
\newtheorem{defn}[thm]{Definition}
\newtheorem{rmk}[thm]{Remark}

\newcommand{\demo}{ {\it   Proof. }}

\begin{document}

\vspace*{2cm}

\subjclass{Primary 20E26, 20F65. Secondary 20E05}

\title[Subgroup conjugacy separability]{Subgroup conjugacy separability\\ for surface groups}

\author{Oleg Bogopolski}
\address{{Sobolev Institute of Mathematics of Siberian Branch of Russian Academy
of Sciences, Novosibirsk, Russia}\newline {and D\"{u}sseldorf University, Germany}}
\email{Oleg$\_$Bogopolski@yahoo.com}

\author{Kai-Uwe Bux}
\address{Bielefeld University, Germany}
\email{bux$\_$2009@kubux.net}


\thanks{$\dagger$ 
This research was partially supported by
SFB 701, ``Spectral Structures and Topological Methods in Mathematics'', at Bielefeld University.}


\begin{abstract}

A group $G$ is called subgroup conjugacy separable (abbreviated as SCS), if
any two finitely generated and non-conjugate subgroups of
$G$ remain non-conjugate in some finite quotient of $G$.
We prove that free groups and the fundamental groups of  orientable closed compact surfaces are SCS.

\end{abstract}

\maketitle

\section{Introduction}

Let $G$ be a group. We call $G$ {\it subgroup conjugacy separable} (SCS) if the following
condition holds:

\medskip
For any two finitely generated subgroups $H_1$ and $H_2$ that are not conjugate in
$G$, there is a homomorphism of $G$ onto a finite quotient $\overline{G}$ such that the images of $H_1$ and $H_2$
are not conjugate in $\overline{G}$.

\medskip

This property logically continues the following series of well
known properties of groups: residual finiteness, conjugacy separability, and subgroup separability (LERF).
Note that SCS-groups are residually finite, but there are residually finite, and even conjugacy separable groups,
which are not SCS-groups. The SCS-property is relatively new and there is not much known about, which groups enjoy this property.

We know of only two papers on SCS: In~\cite{GS}, F.~Grune\-wald and D.~Segal proved that all virtually polycyclic groups are SCS (see also Theorem~7 in Chapter~4 of~\cite{Segal}). In the preprint~\cite{BG}, O. Bogopolski and F.~Grune\-wald proved that free groups and some virtualy free groups are SCS:

\begin{thm}\label{scs_free} {\rm (see~\cite{BG})}
Free groups are SCS.
\end{thm}

The main goal of this paper is to show:

\begin{thm}\label{main_scs}
Fundamental groups of orientable closed compact surfaces are SCS.
\end{thm}

The main tool for establishing these theorems are Theorems~\ref{sics_free} and~\ref{main.sics}. We use there
the following definitions:

For two subgroups $A$ and $B$ of a group $C$, we say that $A$ is {\it conjugate into} $B$ if
$A^c\leqslant B$ for some element $c\in C$. Here $A^c=c^{-1}Ac$.
A group $G$ is called {\it subgroup into conjugacy separable} (SICS) if the following
condition holds:


For any two finitely generated subgroups $H_1$ and $H_2$ such that $H_2$ is not conjugate into $H_1$ in $G$, there is a homomorphism from $G$ onto a finite quotient $\overline{G}$ such that the image of $H_2$ is not conjugate into the image of $H_1$ in $\overline{G}$.

\medskip


\begin{thm}\label{sics_free} {\rm (see~\cite{BG})} Free groups are SICS.
\end{thm}

\begin{thm}\label{main.sics} Fundamental groups of orientable
closed compact surfaces are SICS.
\end{thm}

In Corollary~\ref{reduction}, we reduce Theorems~\ref{scs_free} and~\ref{main_scs} to Theorems~\ref{sics_free} and~\ref{main.sics}. The proof of Theorem~\ref{sics_free} is short and uses simple graphical arguments. Since this proof was never published in a journal and since we use some ideas hidden in it, we decided to include it in this paper for completeness.
To the contrast, the proof of Theorem~\ref{main.sics} is long and uses additional ideas based on the Scott realization theorem, on the Hurwitz realization problem, and on the Bass -- Serre theory of groups acting on trees.
The proof of this theorem will account
for the main part of the paper. We put technical lemmas on covering of surfaces into two appendices.

\section{Chain Condition}

We say that an automorphism $\alpha:G\rightarrow G$ of the group $G$ is {\it expanding} if there is a finitely
generated subgroup $H\leqslant G$ such that $H < H^{\alpha}$ is a strict inclusion.

\medskip




\medskip

Not every group admits expanding automorphisms. E.g., the following result of Takahasi
implies that free groups of finite rank do not.

\begin{thm}\label{Takahasi_free}
{\rm (see \cite{Tak})} Let
$$H_1 < H_2 < H_3 < \dots \leqslant F_n$$
be a strictly ascending infinite chain of finitely generated subgroups of the finitely generated
free group $F_n$.  The sequence of ranks of free groups $H_i$ is unbounded.
\end{thm}

\medskip

We can promote Takahasi's Theorem to surface groups:

\begin{prop}\label{cor_bux_2} Let $S$ be a compact closed surface and let
$$H_1 < H_2 < H_3 < \dots \leqslant \pi_1(S)$$
be a strictly ascending infinite chain of finitely generated subgroups. Then the groups $H_i$ are free
and the sequence of
their ranks is unbounded.
\end{prop}

{\it Proof.} The subgroup $H :=\cup_i H_i$
cannot be finitely generated as the
chain is strictly ascending. Hence $H$ has infinite index in $\pi_1(S)$ and is free of
infinite countable rank. Embedding $H$ into $F_2$, we have the
strictly ascending chain
$$H_1 < H_2 < H_3 < \dots \leqslant F_2,$$
and the claim follows from Theorem~\ref{Takahasi_free}.\hfill $\Box$




\begin{cor}\label{no exp autos} Finitely generated free groups and fundamental groups of compact closed surfaces
do not admit expanding automorphisms.
\end{cor}

\demo  Let $G$ be such a group. If $\alpha : G \rightarrow G$ were expanding on a subgroup $H$,
we would have the strictly ascending chain
$$H < H^{\alpha} < H^{\alpha^2} < \dots $$
All groups in this chain are isomorphic and hence their free ranks coincide, contradicting
Theorem~\ref{Takahasi_free} or Proposition~\ref{cor_bux_2}. \hfill $\Box$


\begin{lem}\label{observation}
Let $H_1$ and $H_2$ be two finitely generated
subgroups of a group $G$.
Assume that $G$ does not admit
expanding inner automorphisms. If $H_2$ conjugates into $H_1$ and $H_1$ conjugates into $H_2$,
then $H_1$ and $H_2$ are conjugate. More precisely, for any two elements $g,h\in G$ with
$H_2^g\leqslant H_1$ and $H_1^h\leqslant H_2$ one already has equality: $H_2^g=H_1$ and $H_1^h=H_2$.

\end{lem}


{\it Proof.} We have $H_1\leqslant H_2^{h^{-1}}\leqslant H_1^{g^{-1}h^{-1}}$.
Put $f := g^{-1}h^{-1}$ and consider the associated
inner automorphism. Since it is not expanding, the inclusion $H_1 \leqslant H_1^f$
is not strict. Hence
$H_1 = H_1^f$ that implies $H_1^h=H_2$ and $H_2^g=H_1$.\hfill $\Box$

\begin{prop}
Suppose that a group $G$ does not admit
expanding inner automorphisms and is SICS, then it is SCS.
\end{prop}

{\it Proof.} Let $H_1$ and $H_2$ be two non-conjugate finitely generated subgroups of $G$. By
Lemma~\ref{observation}, $H_1$ is not conjugate into $H_2$ or $H_2$ is not conjugate into $H_1$. Both cases
are symmetric and we assume that $H_2$ is not conjugate into $H_1$.
Since $G$ is SICS, there exists a homomorphism from $G$ onto a finite group $\overline{G}$
such that the image of $H_2$ is not conjugate into the image of $H_1$ in $\overline{G}$.
In particular, the image of $H_2$ is not conjugate to the image of $H_1$ in $\overline{G}$. Hence, $G$ is SCS.
\hfill $\Box$


\medskip

From this and Corollary~\ref{no exp autos}, we get:

\begin{cor}\label{reduction_free}
If a free group is SICS, then it is SCS.
\end{cor}

\begin{cor}\label{reduction}
If the fundamental group $\pi_1(S)$ of a closed
compact surface $S$ is SICS then it is SCS.\hfill $\Box$
\end{cor}

\section{Auxiliary statements}

For a compact formulation of further results, we need the following terminology:
Let $H_1$ and $H_2$ be finitely generated subgroups of $G$. We say that $H_1$ is
{\it con-sepa\-rated from $H_2$ within $G$} if there is a finite index subgroup $D\leqslant G$ containing $H_1$
such that $H_2$ is not conjugate into $D$. We call $D$ a {\it witness} of con-separation. Note that being
con-separated is not a symmetric relation. We call $H_1$ {\it con-separated in} $G$ if $H_1$
is con-separated from any finitely generated subgroup $H_2\leqslant G$ that is not already conjugate into $H_1$.

\begin{lem}\label{propJuly} Let $H_1,H_2$ be subgroups of a group $G$. Then the following conditions are equivalent:

\begin{enumerate}
\item[{\rm (1)}] $H_1$ is con-separated from $H_2$ within $G$;

\item[{\rm (2)}] There exists a homomorphism from $G$ onto a finite group $\overline{G}$ such that
the image of $H_2$ is not conjugate into the image of $H_1$.
\end{enumerate}
\end{lem}

{\it Proof.} $(1)\Rightarrow (2)$: Let $D$ be the witness of con-separation for $H_1$ from $H_2$.
Then $D$ contains a finite index subgroup $N$ which is normal in $G$.
Obviously, the image of $H_2$ in $G/N$ is not conjugate into the image of $H_1$.

$(2)\Rightarrow (1)$: If $\varphi:G\rightarrow \overline{G}$ is the homomorphism from (2), then $D:=H_1\cdot {\text{\rm ker}}\varphi$ is the witness of con-separation for $H_1$ from $H_2$.
$\hfill$ $\Box$

\medskip


The following lemma enables to push the con-separability within
a finite index subgroup to the con-separability within the whole group.

\begin{lem}\label{overgroups}
Let $G'$ be a finite index subgroup of $G$ and
let $H_1$, $H_2$ be two finitely generated subgroups of $G'$. Let $g_1,\dots ,g_k$
be a set of representatives for the left cosets $G/G'$.
If $H_1$ is con-separated from $H_2^{g_i}$ in $G'$ for each $i$ such that $H_2^{g_i}$
is a subgroup of $G'$, then $H_1$ is con-separated from $H_2$ in $G$.
In particular, if $H_1$ is con-separated in $G'$ it is also con-separated in $G$.
\end{lem}

{\it Proof.} If $H_2^{g_i}$ is contained in $G'$, let $D_i\leqslant G'$ be a witness that $H_1$ is con-separated from
$H_2^{g_i}$ within $G'$. Otherwise, put $D_i := G'$. Note that in either case, $H_2^{g_i}$
is not conjugate into $D_i$ by a conjugating element of $G'$.

\medskip
\noindent We claim that $H_1$ is con-separated from $H_2$ within $G$ with witness
$D := D_1 \cap \dots \cap D_k$. For contradiction, assume $H_2^g\leqslant D$ for some
$g\in G$. We write $g = g_ih$
for some $h\in G'$. Then $H_2^{g_ih}\leqslant D\leqslant D_i$ whence $H_2^{g_i}$
would be conjugate into $D_i$ by a
conjugating element of $G'$. This is a contradiction.\hfill $\Box$

\section{Free groups are SICS}

\subsection{Notations}

Let $\Gamma$ be a graph. By $\Gamma^0$ we denote the set of its vertices and by $\Gamma^1$ the set of its edges.
The inverse of an edge $e\in \Gamma^1$ is denoted by $\overline{e}$, the initial and the terminal vertices of $e$ are denoted by $i(e)$ and $t(e)$. A nontrivial path in $\Gamma$ is a sequence of edges $e_1e_2\dots e_m$, where $m\in \mathbb{N}$ and $t(e_{s})=i(e_{s+1})$ for $s=1,\dots,m-1$. Any vertex of $\Gamma$ is considered as a trivial path.
The inverse to a path $p$ is denoted by $\overline{p}$.

Let $F$ be a free group with finite basis $x_1,\dots,x_n$.
Let $R$ be the graph consisting of one vertex $v$ and $n$ oriented edges $e_1,\dots,e_n$.
We label $e_j$ by $x_j$ and $\overline{e}_j$ by $x_j^{-1}$.
We will identify $F$ with $\pi_1(R,v)$ by identifying $x_j$ with the homotopy class $[e_j]$.

To every subgroup $H\leqslant F$ corresponds a covering map $\varphi: (\Gamma_H,v_H)\rightarrow (R,v)$,
such that $H$ is the image of the induced map $\varphi_{\ast}:\pi_1(\Gamma_H,v_H)\rightarrow \pi_1(R,v)$.
We lift the labeling of $R$ to $\Gamma_H$. So, an edge $e$ of $\Gamma_H$ is labeled  by $x$ if its image $\varphi(e)$ is labeled by $x$.

If $H$ is finitely generated, then $\Gamma_H$ has a finite {\it core}, ${\text{\rm Core}}(\Gamma_H)$, i.e. a finite connected subgraph which is homotopy equivalent to $\Gamma_H$.
We can enlarge ${\text{\rm Core}}(\Gamma_H)$ if necessary and assume that $v_H$
is a vertex of ${\text{\rm Core}}(\Gamma_H)$
and that every vertex of ${\text{\rm Core}}(\Gamma_H)$ has valency 1 or $2n$.
The vertices of valency~1 and the edges of ${\text{\rm Core}}(\Gamma_H)$ emanating from these vertices are called {\it outer}.
For any outer edge $e$, there is a unique reduced path $e_1e_2\dots e_k$
in ${\text{\rm Core}}(\Gamma_H)$ such that $e_1=e$, the labels of the edges $e_j$ coincide and the last edge $e_k$  terminates at an outer vertex.
We will write $e_1^{\,\,\ast}=e_k$.
We will say that the outer edges $e$ and $\overline{e^{\ast}}$ are {\it associated}. This defines an equivalence relation on the set of all outer edges.
From each pair of associated edges we chose one representative and
let $\mathcal{E}$ be the set of all such representatives.

\subsection{Proof of Theorem~\ref{sics_free}}

First we reduce the theorem to the case where the ambient free group is finitely generated.
Let $X$ be an arbitrary set and $F(X)$ be the free group with the basis $X$.
Let $H_1,H_2$ be two finitely generated subgroups of $F(X)$ such that $H_2$ is not conjugate into $H_1$.
Clearly, there exists a finite subset $Y\subseteq X$ such that $\langle H_1,H_2\rangle\leqslant F(Y)$.
If Theorem~\ref{sics_free} holds for finitely generated free groups,
then there exists a finite quotient of $F(Y)$ where the image of $H_2$ is not conjugate into the image of $H_1$.
Since $F(Y)$ is a free factor of $F(X)$, there exists a finite quotient of $F(X)$ with the same property.

So, let $F$ be a free group with finite basis $x_1,\dots,x_n$.
Let $H_1,H_2$ be two finitely generated subgroups of $F$ such that $H_2$ is not conjugate into $H_1$.
By Lemma~\ref{propJuly}, it suffices to construct a finite index subgroup $D$ of $F$ that contains $H_1$ and does not contain a conjugate of $H_2$.

Let $H_2=\langle h_1,\dots, h_r\rangle$, where all $h_j$ are nontrivial, and let $C=\max\{ |h_1|,\dots ,|h_r|\}$,
where $|h|$ denotes the length of $h$ with respect to the given basis of $F$. Since $F$ is residually finite,
there exists a normal subgroup $K$ of finite index in $F$ such that $K$
does not contain nontrivial elements of $F$ of length $C$ or smaller. Since $K$ is normal, $K$ does not contain any conjugate to these elements.
This means that the covering graph $\Gamma_{K}$ is finite, every its vertex has valency $2n$, and

$${\text {\it every cycle in}}\hspace*{2mm} \Gamma_{K}\hspace*{2mm} {\text {\it has length at least}}\hspace*{2mm} C+1.\eqno{(4.1)}
$$


Now we will embed ${\text{\rm Core}}(\Gamma_{H_1})$ into a finite labeled graph $\Delta$ without outer edges.
For every edge $e\in \mathcal{E}$ we choose an edge $\widehat{e}$ in $\Gamma_{K}$ with the same label.
Let $\Delta$ be the labeled graph obtained from the disjoint union of graphs

$${\text{\rm Core}}(\Gamma_{H_1})\bigsqcup\,\, \underset{e\in \mathcal{E}}{\sqcup}\,\, \Bigl(\Gamma_{K}\setminus\{\widehat{e},\,\overline{\widehat{e}}\}\Bigr)\eqno{(4.2)} $$
by identifying the vertices $i(e)$ with $i(\widehat{e})$ and $t(e^{\ast})$
with $t(\widehat{e})$ for every $e\in \mathcal{E}$.
It follows from (4.1) that

\medskip

{\it the distance in $\Gamma_{K}\setminus\{\widehat{e},\,\overline{\widehat{e}}\}$
between the vertices $i(\widehat{e})$ and $t(\widehat{e})$ is at least $C$.}\hfill $(4.3)$

\medskip

\hspace*{25mm}\includegraphics[scale=0.4]{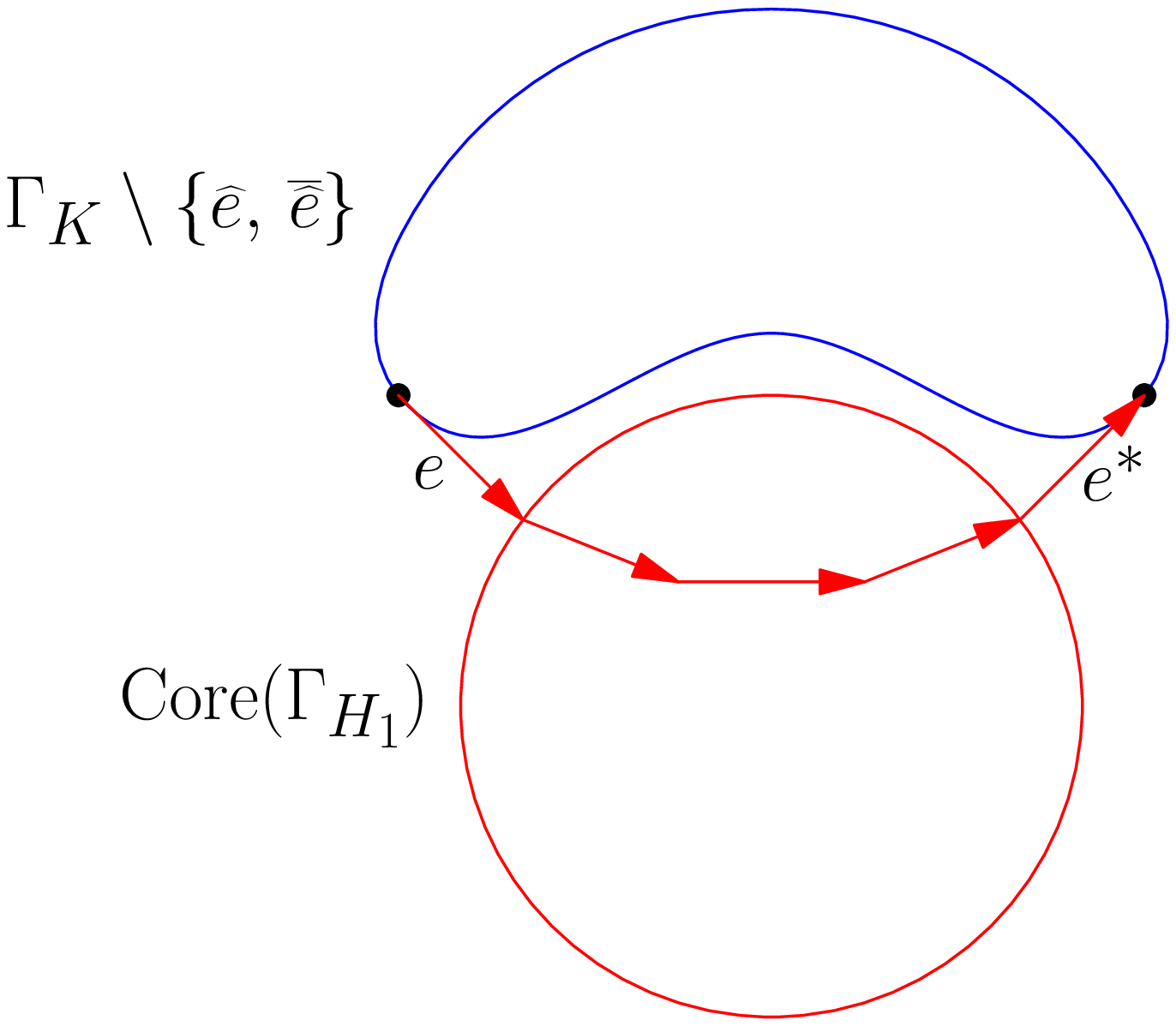}



\vspace*{-65mm}
\begin{center} {\sc Figure 1.} The graph $\Delta$
\end{center}

Since every vertex of $\Delta$ has valency $2n$,
there is a finitely sheeted covering map $\psi: (\Delta, v_{H_1})\rightarrow (R,v)$ respecting the labeling.
Thus $\Delta=\Gamma_D$ for some finite index subgroup $D$ of $F$. Since ${\text{\rm Core}}(\Gamma_{H_1})$ is a subgraph of $\Gamma_D$, the subgroup $D$ contains $H_1$ as a free factor: $D=H_1\ast L$.

We show that $H_2=\langle h_1,\dots ,h_r\rangle$ is not conjugate into $D$.
Suppose the contrary: $gH_2g^{-1}\leqslant D$ for some $g\in F$.
Then the path in $\Delta$ starting at $v_{H_1}$ and having the label $ghg^{-1}$
is closed for each $h\in H_2$.

Let $\gamma$ be the reduced path in $\Delta$ starting at $v_{H_1}$ and having the label $g$, and let $\beta_j$
be the reduced path in $\Delta$ starting at $t(\gamma)$ and having the label $h_j$, $j=1,\dots,r$.
Then each $\beta_j$ is a closed path of length at most $C$. By (4.1) and (4.3), we have the following properties:

\vspace*{-1mm}
\begin{enumerate}
\item[{\small${\bullet}$}]  $\beta_j$ does not lie completely in $\Gamma_{K}\setminus\{\widehat{e},\,\overline{\widehat{e}}\}$ for any $e\in \mathcal{E}$.

\item[{\small${\bullet}$}] $\beta_j$ has no nontrivial subpaths which lie in some $\Gamma_{K}\setminus\{\widehat{e},\,\overline{\widehat{e}}\}$ and have endpoints in
$\{ i(\widehat{e}), t(\widehat{e})\}$.
\end{enumerate}
We consider two cases:

\medskip

{\it Case 1.} Suppose that $t(\gamma)$ lies in some $\Gamma_{K}\setminus\{\widehat{e},\,\overline{\widehat{e}}\}$.
Then there exists a path $\beta$ in $\Gamma_{K}\setminus\{\widehat{e},\,\overline{\widehat{e}}\}$  from $t(\gamma)$ to an outer vertex of ${\text{\rm Core}}(\Gamma_{H_1})$, such that each $\beta_j$ has the form $\beta_j=\beta\beta_j'\bar{\beta}$ with $\beta_j'$ in ${\text{\rm Core}}(\Gamma_{H_1})$.

Let $\delta$ be a path in ${\text{\rm Core}}(\Gamma_{H_1})$ from $v_{H_1}$ to $t(\beta)$.
Then each $\delta\beta'_j\bar{\delta}$ is a closed path in ${\text{\rm Core}}(\Gamma_{H_1})$ based at $v_{H_1}$.
Therefore the label of each path $\delta\bar{\beta}\beta_j\beta\bar{\delta}$ lies in $H_1$.
Let $g_1$ be the label of $\delta\bar{\beta}$. Then $g_1h_jg_1^{-1}\in H_1$ for each $j=1,\dots,r$.
Hence, $g_1H_2g_1^{-1}\leqslant H_1$, a contradiction.

\medskip

{\it Case 2.} Suppose that $t(\gamma)$ lies in ${\text{\rm Core}}(\Gamma_{H_1})$. Then all $\beta_j$ lie in
${\text{\rm Core}}(\Gamma_{H_1})$. Let $\delta$ be a path in ${\text{\rm Core}}(\Gamma_{H_1})$ from $v_{H_1}$ to $t(\gamma)$. Then $\delta\beta_j\bar{\delta}$ is a closed path in ${\text{\rm Core}}(\Gamma_{H_1})$ based at $v_{H_1}$  and we obtain a contradiction as above.

\hfill $\Box$

\medskip

Actually, we have proved the following theorem which is stronger than Theorem~\ref{sics_free}.

\begin{thm}\label{by_the_way}
Let $F$ be a finitely generated free group with a basis $X$. For any finitely generated subgroup $H$ of $F$ and for any constant $C>0$, one can construct a subgroup $D\leqslant F$ of finite index such that:

\begin{enumerate}

\item[1)] $\Gamma_D$ contains ${\text{\rm Core}}(\Gamma_H)$ as a subgraph;

\item[2)] any loop in $\Gamma_D$ of length at most $C$ is freely homotopic to a loop in ${\text{\rm Core}}(\Gamma_H)$;

\item[3)] any path in $\Gamma_D$ of length at most $C$ with endpoints in ${\text{\rm Core}}(\Gamma_H)$ is homotopic (with respect to these points) to a path in ${\text{\rm Core}}(\Gamma_H)$.
\end{enumerate}





\end{thm}

We will use this theorem in Section~9.

\section{Branched coverings of graphs}

The {\it girth} of a graph is the length of a shortest circuit therein without backtracking.
In a forest, such circuits do not exist and we say that forests have infinite girth.


\begin{lem}\label{lemma_5.1} Let $\Gamma$ be a finite connected graph and $m$ be a natural number. Then there exists a finite connected graph $\widehat{\Gamma}$ which covers $\Gamma$ and has the girth larger than~$m$.
\end{lem}

{\it Proof.} Let $F$ be the fundamental group of $\Gamma$ with respect to some point.
We choose a maximal tree $T$ in $\Gamma$.
Then $F$ has the basis $X$, which corresponds to the set of edges of $\Gamma$ outside of $T$.
The length function on $F$ with respect to $X$ will be called $X$-length.

Since $F$ is residually finite,
$F$ has a finite index subgroup $H$, which does not contain nontrivial elements of length up to $m$.
Let $N$ be a finite index normal subgroup in $F$ which is contained in $H$.
Then $N$ does not contain conjugates to nontrivial elements of length up to $m$. Let $\widehat{\Gamma}$ be the covering of $\Gamma$ corresponding to~$N$.

Suppose that there exists a simple closed curve in $\widehat{\Gamma}$, whose edge-length is smaller than or equal to $m$. Then the edge-length of its projection in $\Gamma$ is smaller than or equal to $m$. In particular, $N$ contains a nontrivial element of $X$-length smaller than or equal to $m$, a contradiction.\hfill $\Box$



\begin{defn}
{\rm
A {\it $k$-sheeted branched cover} of a graph $\Gamma$ is a map of graphs $p:\widetilde{\Gamma}\rightarrow \Gamma$
satisfying the following conditions:
\begin{enumerate}
\item[1)] Every open edge of $\Gamma$ is covered by $k$ open edges in $\widetilde{\Gamma}$.
\item[2)] For every vertex $\widetilde{v}\in \widetilde{\Gamma}$ and any two edges $e$ and $e'$ emanating from the vertex $v=p(\widetilde{v})$, the number of edges covering $e$ and emanating from $\widetilde{v}$ equals the number of edges covering $e'$ and emanating from $\widetilde{v}$. We call this number $d(\widetilde{v})$ the {\it branched degree} of $\widetilde{v}$.
\item[3)] For any vertex $v\in \Gamma$ the branched degrees of all its preimages add up to $k$.
\end{enumerate}
}
\end{defn}

If all vertices in a branched
cover $\widetilde{\Gamma}\rightarrow \Gamma$ have branched degree 1, then we have an ordinary covering map of graphs.

\begin{prop}\label{cut vertex} Let $\Gamma$ be a finite connected graph. For any vertex $v\in \Gamma$ fix a natural number $d_v$.
Then the following holds:

\begin{enumerate}
\item[{\rm 1)}] There is a branched cover $p:\widetilde{\Gamma}\rightarrow \Gamma$ such that $\widetilde{\Gamma}$ is a finite
connected graph and any vertex $\widetilde{v}\in \widetilde{\Gamma}$ has the branched degree $d_{p(\widetilde{v})}$.

\vspace*{2mm}

\item[{\rm 2)}] Fix a vertex $v\in \Gamma$. If $d_u\geqslant 2$ for each vertex~$u$ at distance 1 from $v$, then the branched cover in {\rm 1)} can be chosen so that some vertex covering $v$ is not a cut vertex.
\end{enumerate}
\end{prop}

{\it Proof.}
1) Remove from $\Gamma$ the middle points of all edges of $\Gamma$.
We get the set of open $1/2$-neighborhoods of vertices of $\Gamma$, say $\mathcal{O}_{v}$, where $v\in \Gamma^0$.
The natural compactification of $\mathcal{O}_v$ is denoted by $\overline{\mathcal{O}}_v$.
For $k\in \mathbb{N}$, let $k\overline{\mathcal{O}}_v$ denote the wedge of $k$ copies of $\overline{\mathcal{O}}_v$ over $v$, i.e. the union of $k$ copies of $\overline{\mathcal{O}}_v$ followed by the identification of all copies of $v$ into one vertex.
Let $d:=\underset{v\in \Gamma^0}{\prod} d_v$.
Then the desired $\widetilde{\Gamma}$ can be constructed by an appropriate gluing of $d/d_v$ copies of $d_v\overline{\mathcal{O}}_v$, where $v$ runs over $\Gamma^0$, and by passing to a connected component.

2) Passing to a further covering if needed, we may assume that the graph
$\widetilde{\Gamma}$ in 1) has girth at least 3. We claim that $p:\widetilde{\Gamma}\rightarrow \Gamma$ is the desired branched covering.

To the contrary, suppose that every vertex in $\widetilde{\Gamma}$ covering $v$ is a cut vertex in $\widetilde{\Gamma}$.
Let $\widetilde{V}\subset \widetilde{\Gamma}$ be the set of all vertices which cover $v$.

Let $\Delta$ be the graph obtained from $\widetilde{\Gamma}$ by replacing each component of $\widetilde{\Gamma}\setminus \widetilde{V}$ by a vertex and each vertex $\widetilde{v}\in \widetilde{V}$ by an edge. Since each $\widetilde{v}\in \widetilde{V}$ is a cut vertex of $\widetilde{\Gamma}$,
the graph $\Delta$ is a finite tree.
Hence $\Delta$ contains a vertex of valency 1.

Therefore, there is a component $C$ of $\widetilde{\Gamma}\setminus \widetilde{V}$ whose closure in $\widetilde{\Gamma}$ contains only one vertex from $\widetilde{V}$, say $\widetilde{v}$.
Let $\widetilde{u}$ be a vertex in $C$ at distance 1 from $\widetilde{v}$. By assumption, we have $d(\widetilde{u})\geqslant 2$. Hence there are at least two vertices at distance 1 from $\widetilde{u}$ which cover $v$ (use that the girth of $\widetilde{\Gamma}$
is at least~3). One of them is $\widetilde{v}$, another is contained in $C$, a contradiction. \hfill $\Box$


\section{Graphs of surfaces}

Let $S$ be a compact surface decomposed into subsurfaces $B_1,\dots, B_m$, i.e., $S$ is obtained from the subsurfaces
by gluing them along boundary circles.

This can be regarded as a {\it graph of spaces} decomposition of $S$: we have one vertex for each region $B_i$
and an edge for each identification of a pair of boundary circles. In the surface, these circles appear as {\it cutting circles} along which $S$ is decomposed.

Let $\Gamma$ be underlying graph of the decomposition. There is a continuous projection
$f:S\twoheadrightarrow \Gamma$, where the preimages of edges are annuli (thickenings of the cutting circles)
and preimages of vertices are regions $B_i$ (up to a missing collar around the boundary circles).

\begin{construction} {\bf(Pullback along a covering).} {\rm Let $f:S\rightarrow \Gamma$ be as above and let $p:\widetilde{\Gamma}\rightarrow \Gamma$ be a graph covering. Then
$$\widetilde{S}:=\{(s,\widetilde{g})\in S\times \widetilde{\Gamma}\,|\, f(s)=p(\widetilde{g})\}$$
is a surface and the projection onto the first coordinate, $\widetilde{S} \rightarrow S$, is a surface covering. The projection onto the second coordinate is a map $\widetilde{S}\rightarrow \widetilde{\Gamma}$ that corresponds to a graph of space decomposition of the surface $\widetilde{S}$. Note that the two projections make the diagram

$$
\begin{CD}
\widetilde{S} @>>> \widetilde{\Gamma}\\
@VVV @VVV\\
S @>>> \Gamma
\end{CD}\eqno{(6.1)}
$$
\vspace*{2mm}
commutative.
}
\end{construction}

We can extend the pull-back construction to branched covers, provided we have pieces to cover the vertex spaces
(see below).

\begin{defn}
{\rm Let $B$ be a compact surface with boundary. We call a $d$-sheeted covering map $\widetilde{B}\rightarrow B$ {\it regular} if $\widetilde{B}$ is connected and each boundary circle of $\widetilde{B}$ is mapped homeomorphically to the corresponding boundary circle of $B$. Equivalently, any boundary circle of $B$ is covered by precisely $d$
boundary circles of $\widetilde{B}$.
}
\end{defn}

Note that the identity map $B\rightarrow B$ is a regular 1-sheeted covering.

\begin{construction} {\bf (Regular coverings)}
{\rm Let $S$ be a compact surface with boundary. By capping off the boundary circles with discs, we obtain a closed surface $\overline{S}$. We call the discs attached to the boundary circles {\it caps}. Suppose that $\widetilde{\overline{S}}\rightarrow \overline{S}$ is a $d$-sheeted covering.
Since a cap is simply connected, it has $d$ preimages all of which are discs in $\widetilde{\overline{S}}$.
Removing these preimage discs, we obtain a subsurface $\widetilde{S}$ and a regular $d$-sheeted covering $\widetilde{S}\rightarrow S$.
}
\end{construction}

\begin{rmk}
{\rm The capped surface $\overline{S}$ does allow a $d$-sheeted covering for any $d\geqslant 1$ unless $\overline{S}$
is a sphere or a projected plane.
}
\end{rmk}

\begin{rmk}
{\rm Let $B$ be a compact surface (possibly with boundary) and let $\overline{B}$ be the surface obtained by capping off the boundary components, i.e., attaching a disc to each of the boundary components. Let
$\{\gamma_1,\dots ,\gamma_r\}$ be a finite collection of non-trivial loops in $B$ each of which stays non-trivial in the capped off surface $\overline{B}$. Since surface groups are residually finite, there is a finite sheeted normal covering $\widetilde{\overline{B}}\rightarrow \overline{B}$ such that no loop $\gamma_i$ has a closed lift in $\widetilde{\overline{B}}$. Hence there is a regular covering $\widetilde{B}\rightarrow B$ such that no loop $\gamma_i$ has a closed lift in $\widetilde{B}$.

}
\end{rmk}

\begin{defn}
{\rm Let $f:S\rightarrow \Gamma$ be a graph of spaces decomposition of the surface $S$.
For any vertex $v\in \Gamma$ let $B_v$ denote the corresponding subsurface of $S$.
Let $p:\widetilde{\Gamma}\rightarrow \Gamma$ be a $k$-sheeted branched cover of $\Gamma$.
A covering map $h:\widetilde{S}\rightarrow S$ is called {\it compatible with} $p:\widetilde{\Gamma}\rightarrow \Gamma$ if there exists a graph of spaces decomposition
$\widetilde{f}:\widetilde{S}\rightarrow \widetilde{\Gamma}$ such that the following holds:
\begin{enumerate}
\item[1.] The diagram~(6.1) commutes.
\item[2.]  For any vertex $\widetilde{v}\in \widetilde{\Gamma}$, the induced map $\widetilde{B}_{\widetilde{v}}\rightarrow B_{p(\widetilde{v})}$ is a $d(\widetilde{v})$-sheeted regular covering map. (Here $\widetilde{B}_{\widetilde{v}}$ is the vertex space of $\widetilde{S}$ which lies over the vertex $\widetilde{v}$ of $\widetilde{\Gamma}$.)
\end{enumerate}
}
\end{defn}

\begin{construction}\label{constr} {\bf (Blowing up vertex spaces)} {\rm Let again $f:S\rightarrow \Gamma$ be a graph of spaces
decomposition of the surface $S$. We assume that for each vertex space $B_v$ the capped off surface $\overline{B}_v$ is different from a sphere or a projective plane.
Let $p:\widetilde{\Gamma}\rightarrow \Gamma$ be a $k$-sheeted branched cover of $\Gamma$.
Our goal is to construct a covering $\widetilde{S}\rightarrow S$ compatible with $p:\widetilde{\Gamma}\rightarrow \Gamma$.

For any vertex $\widetilde{v}\in \widetilde{\Gamma}$ let $q_{\widetilde{v}}:\widetilde{B}_{\widetilde{v}}\rightarrow B_{p(\widetilde{v})}$ be a regular $d(\widetilde{v})$-sheeted covering map, where $d(\widetilde{v})$ is the branched degree of $\widetilde{v}$. Then we obtain $\widetilde{S}$ by gluing
the pieces $\widetilde{B}_{\widetilde{v}}$ along boundary circles as the edges in $\widetilde{\Gamma}$ dictate.
Note that there may be several gluing schemata.
}
\end{construction}

\section{From surface groups to topology}

Our goal is to show that the fundamental groups of orientable closed compact surfaces are SCS. By Corollary~\ref{reduction},
it suffices to show that such groups are SICS.
In view of Lemma~\ref{propJuly}, it suffices to show the following:

\begin{thm}\label{oleg1}
Let $S$ be an orientable closed compact surface with basepoint $\ast$. Any finitely generated subgroup $H_1$ of
the fundamental group $\pi_1(S,\ast)$ is con-separated.
\end{thm}

Here, we reduce Theorem~\ref{oleg1} to a purely
topological statement (Problem~\ref{problem}), which will be dealt in the next section.


Let $S$ be a closed surface. A subgroup $H$ of its fundamental group is called {\it geometric} if there is
an incompressible subsurface $A\subseteq S$ containing $\ast$, whose fundamental group with respect to $\ast$ is $H$.
(Recall that a compact subsurface $A$ of the surface $S$ is called {\it incompressible} if the embedding $A\hookrightarrow S$ induces the embedding of fundamental groups: $\pi_1(A)\hookrightarrow \pi_1(S)$. Geometrically this means that the complement of $A$ in $S$ does not have disc components.)

It suffices to show that geometric subgroups of $\pi_1(S,\ast)$ are con-separated. Indeed, by~\cite[Theorem~3.3]{Scott},
there is a finite sheeted cover $\widetilde{S}$ of $S$ such that $H_1$ is a geometric subgroup of $\pi_1(\widetilde{S},\widetilde{\ast})$. If $H_1$ is con-separated in $\pi_1(\widetilde{S},\widetilde{\ast})$, then
it is con-sepa\-rated in $\pi_1(S,\ast)$ by Lemma~\ref{overgroups}.

Hence, we assume from now on without loss of generality that $H_1$ is geometric in $S$.
We also assume that $S$ is a hyperbolic surface, otherwise the theorem is obvious.
Let $A\subset S$ be the subsurface realizing $H_1$.

The subsurface $A$ together with its complementary regions induce a graph of spaces decomposition of $S$ with underlying graph $\Gamma$.
Passing to a cover, we may assume that each complementary region {\it has genus},
i.e. it is not a sphere with holes.

Moreover, we can arrange that $A$ has a unique complementary component in $S$. Indeed,
let $v$ be the vertex of $\Gamma$ corresponding to $A$.
By Proposition~\ref{cut vertex}, we can construct a finitely-sheeted branched cover $p:\widetilde{\Gamma}\rightarrow \Gamma$ such that all vertices of $\widetilde{\Gamma}$ covering $v$ have branched degree 1, all other vertices of $\widetilde\Gamma$
have branched degree 2, and some vertex $\widetilde{v}$ over $v$ is not a cut vertex of $\widetilde{\Gamma}$.
Let $\widetilde{S}\rightarrow S$ be a covering map compatible with $p$, see Construction~\ref{constr} (we can use it
since the vertex spaces of $S$ different from $A$ have genus).
Let $\widetilde{A}$ be the subspace of $\widetilde{S}$ corresponding to $\widetilde{v}$.
Then $\widetilde{A}$ has only one complementary region and $\widetilde{A}$ realizes  $H_1$.


\begin{defn} {\rm We say that subsurface $A$ of a surface $T$ has {\it good shape} if $A$
is incompressible, $A$ has only one complementary component~$B$, and $B$ has genus.}
\end{defn}

We have argued:

\begin{lem} Let $A$ be an incompressible subsurface of $S$ realizing $H_1$. Then there is a finitely-sheeted covering $\widetilde{S}\rightarrow S$ that contains a lift of $A$ of good shape.
\end{lem}

When $A$ has good shape, we have a further reduction:

\begin{thm} Let $S$ be an orientable closed compact surface with $\chi(S)\leqslant -1$ and with basepoint $\ast$, and let $H_1$ be a geometric subgroup of $\pi_1(S,\ast)$. Suppose that a subsurface $A\subseteq S$ which realizes $H_1$ has good shape. Then $H_1$ is con-separated, provided the following holds:

\medskip

\hspace*{5mm}\begin{minipage}[r]{120mm} {\sl For any element $g\in \pi_1(S,\ast)$ not conjugate into $H_1$, there exists a
finite index subgroup $D\leqslant \pi_1(S,\ast)$ containing $H_1$ but not containing any conjugate of $g$.}
\end{minipage}

\end{thm}

\noindent
(In this case, we say that $H_1$ is {\it con-separated} from each element $g$ not conjugate into $H_1$,
and we say that $D$ is a {\it witness} for the separation.)

\medskip

\demo Let $H_2\leqslant \pi_1(S,\ast)$ be a finitely generated subgroup not conjugate into $H_1$. If there is an element $g\in H_2$ that cannot be conjugate into $H_1$, any subgroup $D\leqslant \pi_1(S,\ast)$
con-separating $H_1$ from $g$ will also con-separate $H_1$ from $H_2$.
Therefore, it suffices to prove the following claim:

{\it Claim.}
Suppose that each element of $H_2$ is conjugate into $H_1$. Then $H_2$ is conjugate into $H_1$.

\medskip

We may assume that $H_2$ is not cyclic, otherwise the claim is trivial. The subsurface $A$, the closure of its complement $B$, and the common boundary components of $A$ and $B$ induce a graph of spaces decomposition of $S$.
The group $G:=\pi_1(S,\ast)$ acts on the associated Bass-Serre tree $T$. Every element of $H_2$ has a conjugate in $H_1=\pi_1(A,\ast)$ and, hence, acts elliptically. Therefore $H_2$ has a global fixed vertex in~$T$.
If $H_2$ fixes a vertex corresponding to a conjugate of $\pi_1(A)$, then $H_2$ is conjugate to $H_1$ and we are done.

Now suppose that $H_2$ fixes a vertex corresponding to a conjugate of $\pi_1(B)$.
By assumption of the claim, each element $x\in H_2$ fixes a vertex corresponding to a conjugate of $\pi_1(A)$.
Therefore $x$ fixes an edge of $T$. The edges of $T$ correspond to the boundary components $R_1,\dots,R_k$ of $A$.
Let $a_j$ be a generator of $\pi_1(R_j)$, $j=1,\dots ,k$. Then each $x\in H_2$ is conjugate to a power of some $a_j$.

Since $H_2$ is noncyclic, the chain of the commutator subgroups of $H_2$ is strictly descending and infinite:
$H_2=H_2^{(0)}>H_2^{(1)}>H_2^{(2)}>\dots $, and $\overset{\infty}{\underset{i=0}{\cap}} H_2^{(i)}=1$.
For any natural $n$, we choose a nontrivial $x_n\in H_2^{(n)}$. Then $x_n\in G^{(f(n))}\setminus G^{(f(n)+1)}$ for some $f(n)\geqslant n$. There is an infinite subset $I\subseteq \mathbb{N}$ such that each $x_i$, $i\in I$, is conjugate to a nonzero power of the same $a_j$, say $x_i\sim a_j^{\ell(i)}$ with $\ell(i)\neq 0$. Then  $a_j^{\ell(i)}\in G^{(f(i))}\setminus G^{(f(i)+1)}$.
Since the quotient $G^{(f(i))}/ G^{(f(i)+1)}$ is torsionfree, we have $a_j\in G^{(f(i))}\setminus G^{(f(i)+1)}$ for each $i\in I$, that is impossible since $f(I)$ is infinite.
\hfill $\Box$

\medskip

The remaining problem has a straightforward topological interpretation.

\begin{problem}\label{problem}  {\rm Let $S$ be an orientable closed compact surface with $\chi(S)\leqslant -1$.
Given a subsurface $A\subset S$ of good shape and given a loop $\gamma\subset S$ that cannot be freely homotoped into $A$, find a finitely-sheeted covering $\widetilde{S}\rightarrow S$ such that $A$ lifts but $\gamma$ does not.}
\end{problem}

\section{Proof of Theorem~\ref{main.sics}: Solution of Problem~\ref{problem}
}

We explain here how we construct such $\widetilde{S}$.
Let $R_1,\dots,R_n$ be all boundary components of $A$. We fix a hyperbolic metric on $S$,
such that these boundary components became geodesics (see~\cite{Farb}).
We may assume that $\gamma$ is the shortest curve in the free homotopy class of $\gamma$.
Let $C$ be the length of $\gamma$ with respect to this metric.
Every covering of $S$ inherits the metric on $S$. A curve in a covering of $S$
is called {\it $C$-short}, if its length does not exceed $C$;
otherwise a curve is called {\it $C$-long}. We set $B:={\text{\rm cl}}(S\setminus A)$.

We will construct special coverings $B_1,\dots ,B_{n+2}$ of $B$
and special coverings $A_{n+1}$ and $A_{n+2}$ of $A$,
and then we will construct $\widetilde{S}$ by gluing several copies of these coverings and several copies of $A$
according to the schema described below.

\subsection{Special coverings of $A$ and $B$}

{\rm Let $\bar{A}$ be a covering of $A$. A boundary component of $\bar{A}$
is called an $(R_i,d)$-{\it boundary} if it covers the boundary component $R_i$ of $A$ with degree $d$.
We use the same definition for $B$ instead of $A$.
Sometimes we will shorten this wording to {\it $d$-boundary}.
Note that the surfaces $A$ and $B$ have only 1-boundaries.}

\medskip




By Lemma~\ref{very technical}, there exist coverings $B_1,\dots ,B_{n+2}$ of $B$
and coverings $A_{n+1},A_{n+2}$ of $A$ which satisfy the following two conditions:

\begin{enumerate}
\item[\bf 1.] {\bf Condition on boundaries}:\\
\vspace*{1mm}
There are natural numbers $M,N,N',N''$ such that $M>1$ and the following holds:

\begin{itemize}
\item[i)] For each $i=1,\dots,n$, the surface $B_i$ contains exactly

$M$ boundary components of type $(R_i,1)$,

$N$ boundary components of type $(R_i,M)$,

$N'$ boundary components of type $(R_i,2M)$,

$N''$ boundary components of type $(R_j,2M)$ for each $j\in \{1,\dots,n\}\setminus \{i\}$.

\vspace*{1mm}



\item[ii)] The surfaces $A_{n+1}$ and $B_{n+1}$ contain only
boundary components of type $(R_j,M)$ for each $j=1,\dots,n$.
Clearly, the number of such components for each $j$ is the same.

\vspace*{1mm}

\item[iii)] The surfaces $A_{n+2}$ and $B_{n+2}$ contain only
boundary components of type $(R_j,2M)$ for each $j=1,\dots,n$.
Clearly, the number of such components for each $j$ is the same.

\end{itemize}

\vspace*{2mm}

\item[\bf 2.] {\bf Condition on short curves}:\\ For every $K\in \{A_{n+1},B_{n+1},A_{n+2},B_{n+2} \}\cup\{B_1,\dots,B_{n}\}$ and every $C$-short geodesic curve $\beta$ in $K$ we have:

\begin{itemize}




\vspace*{1mm}

\item[i)] If $\beta$ is a loop, then $\beta$ lies in some $1$-boundary of $K$.
In particular, $A_{n+1},B_{n+1},A_{n+2}$ and $B_{n+2}$ don't contain $C$-short geodesic loops.

\vspace*{1mm}

\item[ii)] If $K\in \{B_1,\dots ,B_n\}$ and the endpoints of $\beta$ lie in 1-boundaries of $K$, then
$\beta$ lies completely in one 1-boundary component of $K$.

\vspace*{1mm}

\item[iii)] If $K\in \{A_{n+1},B_{n+1},A_{n+2},B_{n+2}\}$
and the endpoints of $\beta$ lie on a boundary component of $K$,
then $\beta$ lies completely in this boundary component.
\end{itemize}
\end{enumerate}

\subsection{Construction of $\widetilde{S}$}

{\rm Let $d$ be the minimum of distances between different boundary components of each of the surfaces $A_{n+1},A_{n+2},B_{n+1},B_{n+2}$, and let $T$ be the minimal odd number such that $T> C/d$, where $C$ is the length of $\gamma$.}

The construction up to Step~3 is illustrated by Figure~2.

\begin{itemize}
\item[Step 1.] Recall that, for every $i=1,\dots ,n$, the space $A$ contains a single $(R_i,1)$-boundary
and the space $B_i$ contains $M$ copies of $(R_i,1)$-boundaries.\\
We glue $M$ copies of $A$ to $B_1,\dots, B_n$ along the corresponding $1$-boundaries so that in the resulting
space $S_1$ each copy of $A$ is completely surrounded by $B_1,\dots,B_n$ and $S_1$ has no $(R_i,1)$-boundaries.

Then $S_1$ contains only $M$-boundaries and $2M$-boundaries. These boundaries lie in the $B_i$-subspaces.
Moreover, the number of $(R_i,M)$-boundaries of $S_1$ does not depend on $i$.
The same is true for the $(R_i,2M)$-boundaries of $S_1$.







\item[Step 2.]
To each boundary component of $S_1$ we glue a copy of $A_{n+1}$ or $A_{n+2}$ along the corresponding  boundary component. The resulting space $S_2$ contains only $M$-boundaries and $2M$-boundaries, and these boundaries lie in the $A_{n+1}$-subspaces and in the $A_{n+2}$-subspaces, respectively.







\item[Step 3.]
To each boundary component of $S_2$ we glue a copy of $B_{n+1}$ or $B_{n+2}$ along the corresponding  boundary component. The resulting space $S_3$ contains only $M$-boundaries and $2M$-boundaries, and these boundaries lie in the $B_{n+1}$-subspaces and in the $B_{n+2}$-subspaces, respectively.



\item[Steps 4 to $T$.] We continue the process by applying the procedure described in Steps~2 and~3 alternately to current spaces $S_i$ until we get the space $S_T$, where $T$ is the constant defined above.

    The space $S_T$ contains only $M$- and $2M$-boundaries, and since $T$ is odd,
    these boundaries lie in $B_{n+1}$-subspaces or in $B_{n+2}$-subspaces.
    Moreover, the number of $(R_i,M)$-boundaries of $S_T$ does not depend on $i$,
    and the number of $(R_i,2M)$-boundaries of $S_T$ does not depend on $i$.


\item[Step $T+1$.]  Recall that the number of $(R_i,M)$-boundaries of $A_{n+1}$ does not depend on~$i$,
and the number of $(R_i,2M)$-boundaries of $A_{n+2}$ does not depend on $i$.
We glue several copies of $S_{T}$ and several copies of $A_{n+1}$ and $A_{n+2}$ so that

\begin{enumerate}
\item[a)] the underlying graph is connected and has girth at least $T+1$\\ (use Lemma~\ref{lemma_5.1}),

\item[b)] the resulting space is a closed surface.
\end{enumerate}
\noindent
We denote this surface by $\widetilde{S}$.
\end{itemize}

\vspace*{5mm}


\vspace*{-35cm}
\hspace*{-5cm}
\hspace*{-155mm}\includegraphics[scale=4]{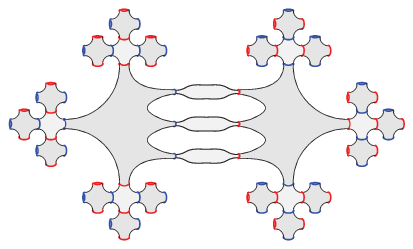}

\vspace*{-735mm}
\begin{center} \hspace*{-10mm}{\sc Figure 2.} An example of the surface $S_3$ with parameters $n=2$, $N=2$, $M=3$.\\ \hspace*{-23mm}For simplicity, we assume here that $N'=N''=0$.
\end{center}

\vspace*{5mm}

 Clearly, $\widetilde{S}$ is a finite cover of $S$ containing a copy of $A$. We shall show that $\gamma$ has no closed lifts in $\widetilde{S}$.

\medskip

{\it Definition.} Let $c$ be a curve in $\widetilde{S}$. A subcurve $c_1$ of $c$ is called a {\it peak} of $c$
if $c_1$ lies in a copy of $K\in \{A_{n+1},B_{n+1},A_{n+2},B_{n+2}\}$
and has endpoints on a boundary component of $K$, and at least one point in the interior of $K$.


\begin{prop}
The geodesic loop $\gamma$ in $S$ has no closed lifts in $\widetilde{S}$.
\end{prop}

{\it Proof.} Suppose that $\widetilde{\gamma}$ is such a lift. Then the length of $\widetilde{\gamma}$ is equal to the length of $\gamma$ which is $C$. Moreover, $\widetilde{\gamma}$ is geodesic as $\gamma$. The curve $\widetilde{\gamma}$ cannot cross different copies of $S_1$, otherwise it crosses at least $(2T-1)$ copies of the spaces $A_{n+1},B_{n+1}, A_{n+2},B_{n+2}$, and the length of $\widetilde{\gamma}$ would be at least $(2T-1)d>C$, a contradiction.
By analyzing the construction of $\widetilde{S}$, we conclude that either $\widetilde{\gamma}$ lies completely in a copy
of $K\in \{A_{n+1},B_{n+1},A_{n+2},B_{n+2}\}$, or
$\widetilde{\gamma}$ has a peak in a copy of $K$, or $\widetilde{\gamma}$ lies completely in a copy of $S_1$. The first is impossible by Condition~2.i), the second by Condition~2.iii).

Thus, we may assume that $\widetilde{\gamma}$ lies in $S_1$. Then $\widetilde{\gamma}$ meets a copy of $A$, otherwise it lies completely in a copy of $B_i$ for some $i\in \{1,\dots,n\}$ that contradicts Condition~2.i).

Suppose that $\widetilde{\gamma}$ intersects the interior of $B_i$.
Then there is a subcurve of $\widetilde{\gamma}$
which intersects the interior of $B_i$ and has endpoints on $1$-boundary components of $B_i$. This contradicts Condition~2.ii).

Hence $\widetilde{\gamma}$ lies completely in a copy of $A$.
But this contradicts the assumption that $\gamma$ cannot be freely homotoped into~$A$.
This shows that $\gamma$ has no closed lifts in $\widetilde{S}$.
\hfill $\Box$


\medskip

Thus, Problem~\ref{problem} is solved. By Section 7, this proves Theorem~\ref{main.sics}.
This and Corollary~\ref{reduction} imply Theorem~\ref{main_scs}.



\section{Appendix A: Constructing metrically large coverings}

\begin{rmk}\label{qi}{\rm
  Let $S$ be a compact nonclosed surface with
  $\chi(S)\leqslant -1$ and with a fixed hyperbolic metric and geodesic
  boundaries. The universal cover $\widehat{S}$ is a convex subset of
  the hyperbolic plane bounded by geodesic lines. It looks like a thickened
  tree.

  The fundamental group $\pi_1(S)$ is free and hence can be identified with the fundamental
  group of a bouquet of circles $\Gamma$. It acts simultaneously on the
  universal covers $\widehat{\Gamma}$ and $\widehat{S}$ by deck
  transformations. There is a $\pi_1(S)$-equivariant homotopy equivalence
  \(
    \hat{f} : \widehat{S} \rightarrow \widehat{\Gamma}
    \text{;}
  \)
  and one may choose
  $\hat{f}$ to be an immersion on boundary components and a
  quasi-isometry, i.e., there is a constant $a\geqslant 1$ such that
  \[
    \frac{1}{a}\cdot d_{\widehat{\Gamma}}(\hat{f}(x),\hat{f}(y))-a
    \leqslant d_{\widehat{S}}(x,y)
    \leqslant a\cdot d_{\widehat{\Gamma}}(\hat{f}(x),\hat{f}(y))+a
  \]
  for any two points $x,y\in\widehat{S}$. Here, $d_{\widehat{S}}$ and
  $d_{\widehat{\Gamma}}$ denote the metrics on $\widehat{S}$ and
  on $\widehat{\Gamma}$.

  For any subgroup $H\leqslant \pi_1(S)$, the homotopy equivalence
  \(
    \widehat{S}/H \rightarrow \widehat{\Gamma} / H
  \)
  induced by $\hat{f}$ is a quasi-isometry for the same
  quasi-isometry constant $a$. It also restricts to immersions on
  boundary components. (In case $H=\pi_1(S)$, this applies to the
  induced map
  \(
    f : S \rightarrow \Gamma
  \).)

  Hence,
  for any covering $\phi:\widetilde{\Gamma}\rightarrow \Gamma$,
  there exists a covering $\psi: \widetilde{S}\rightarrow S$ and a
  homotopy equivalence
  $\widetilde{f}:\widetilde{S}\rightarrow \widetilde{\Gamma}$
  such that $\widetilde{f}$ is a quasi-isometry with quasi-isometry
  constant $a$, and the following diagram is commutative:
  \medskip
  \[
    \begin{CD}
     & \psi \vspace*{-5mm}\\
    \widetilde{S} @ >>> S\\
    \widetilde{f} @VVV @VVV f\\
    \widetilde{\Gamma} @>>> \Gamma\vspace*{-10mm}\\
     & \phi\\
    \end{CD}\eqno{(9.1)}
  \]
  \vspace*{5mm}
}\end{rmk}

Recall that, for any covering $\theta:\widehat{S}\rightarrow S$, a boundary component $\widehat{R}$ of $\widehat{S}$
is called a {\it $k$-boundary} if the restriction of $\theta$ to $\widehat{R}$ has degree $k$. 

{\begin{lem}\label{lemma7} Let $S$ be a compact nonclosed surface with $\chi(S)\leqslant -1$
and with geodesic boundaries with respect to a fixed hyperbolic metric. Let $C_0$ be the maximum of lengths
of boundary components of $S$. For any constant $C>C_0$ and any boundary component $R$ of $S$,
there exist a finite covering $\psi:\widetilde{S}\rightarrow S$ and a boundary component $\widetilde{R}$ of $\widetilde{S}$ such that the following holds.
\begin{enumerate}
\item[{\rm 1)}] $\widetilde{R}$ covers $R$ with degree 1; other boundary components of $\widetilde{S}$ cover the corresponding boundary components of $S$ with degrees larger than 1.
\item[{\rm 2)}] Every $C$-short loop in $\widetilde{S}$ is freely homotopic into $\widetilde{R}$.
\item[{\rm 3)}] Every $C$-short geodesic curve in $\widetilde{S}$ with endpoints on $\widetilde{R}$ lies in $\widetilde{R}$.

\item[{\rm 4)}] For every covering $\theta:\widehat{S}\rightarrow S$, which factors through $\psi:\widetilde{S}\rightarrow S$, the following is satisfied:

    \begin{enumerate}
    \item[(a)] each $C$-short loop in $\widehat{S}$ can be freely homotoped into a boundary component of $\widehat{S}$.
    \item[(b)] every $C$-short geodesic curve in $\widetilde{S}$ with endpoints on a 1-boundary of $\widehat{S}$ lies in this 1-boundary;
    \item[(c)] the distance between any two 1-boundary components of $\widehat{S}$ is larger\break than $C$.
    \end{enumerate}

\end{enumerate}

\end{lem}

{\it Proof.} We are in the situation of the preceding remark and we use its notation.
In particular, $f:S\rightarrow\Gamma$ denotes a homotopy equivalence that restricts to
immersions on boundary components and lifts to a quasi-isometry with constant $a$
on each connected covering space. In addition, we denote by $v$ the basepoint of
the graph $\Gamma$. For any subgroup
\(
  G\leqslant \pi_1(\Gamma,v)
\), let $\Gamma_G$ denote the covering space of $\Gamma$ corresponding to $G$.

By Theorem~\ref{by_the_way}, for any finitely generated subgroup $H$ of $F$ and for any constant $C_1>0$, one can construct a subgroup $D\leqslant F$ of finite index such that:

\begin{enumerate}

\item[i)] $\Gamma_D$ contains ${\text{\rm Core}}(\Gamma_H)$ as a subgraph;

\item[ii)] any loop in $\Gamma_D$ of length at most $C_1$ is freely homotopic to a loop in ${\text{\rm Core}}(\Gamma_H)$;

\item[iii)] any path in $\Gamma_D$ of length at most $C_1$ with endpoints in ${\text{\rm Core}}(\Gamma_H)$ is homotopic
(with respect to these points) to a path in ${\text{\rm Core}}(\Gamma_H)$.
\end{enumerate}


We apply Theorem~\ref{by_the_way} to $F:=\pi_1(S,x)$ and $H:=\pi_1(R,x)$, where $x$ is a point on~$R$, and to $C_1:=aC+a$, and construct the subgroup $D$ of $F$ as above. Let $\phi:(\Gamma_D,v_D)\rightarrow (\Gamma,v)$ be the covering
corresponding to $D$. By Remark~\ref{qi}, there exists
a covering $\psi:(\widetilde{S},\widetilde{x})\rightarrow (S,x)$ such that the right square in the diagram (9.2) is commutative for an appropriate homotopy equivalence $\widetilde{f}$. The commutativity of this square implies that the restriction of $\widetilde{f}$ to any boundary component of $\widetilde{S}$ is an immersion.

Note that the covering $\psi:(\widetilde{S},\widetilde{x})\rightarrow (S,x)$ also corresponds to the subgroup~$D$.
Since $H\leqslant D$, the boundary component $R$ of $S$ has a lift $\widetilde{R}$ covering $R$ with degree~1.
Moreover, $\widetilde{f}$ maps $\widetilde{R}$ to the loop in ${\text{\rm Core}}(\Gamma_H)$.
Thus, we have the following commutative diagram, where $j$ and $i$ are identity embeddings and $\widetilde{f}_{|\widetilde{R}}$ is a homeomorphism:

$$
\begin{CD}
 & j & & \psi \vspace*{-5mm}\\
\widetilde{R} @ >>> \widetilde{S} @ >>> S \\
\widetilde{f}_{|\widetilde{R}} @ VVV @ VVV \hspace*{-16mm}\widetilde{f}\hspace*{14mm} @ VVV f \\
{\text{\rm Core}}(\Gamma_H) @>>> \Gamma_D @>>> \Gamma\vspace*{-10mm}\\
& i  & & \phi\\
\end{CD}
\eqno{(9.2)}
$$

\bigskip

Now one can easily deduce the statements 2) and 3) of the lemma from the statements ii) and iii).
We prove the remaining part of the statement 1). 
Let $\widetilde{R}_1$ be a boundary component of $\widetilde{S}$ different from $\widetilde{R}$. 
By~2), the length of $\widetilde{R}_1$ is larger than $C$. 
Since $C>C_0$, the loop $\widetilde{R}_1$ covers the corresponding $R_1$ with degree larger than 1.
The statement~4) follows from the statements 1)-3).\hfill $\Box$







\begin{lem}\label{7a}
Let $S$ be a compact nonclosed surface with $\chi(S)\leqslant -1$
and with a fixed hyperbolic metric and geodesic boundaries.
For any constant $C>0$, there exists a finite covering $\varphi:\widetilde{S}\rightarrow S$, such that

{\rm 1)} every $C$-short loop in $\widetilde{S}$ is trivial,

{\rm 2)} every $C$-short curve in $\widetilde{S}$ with endpoints on a boundary component of  $\widetilde{S}$  is homotopic, relative to the endpoints, to a geodesic curve in this component.

\end{lem}

{\it Proof.} For any curve $\gamma$ in $S$, let $l(\gamma)$ denote its length.
There are only finitely many closed curves in $S$ (up to free homotopy) of length up to $C$. Furthermore, the fundamental group of $S$ is residualy finite. Hence we can construct a finite covering of $S$ where every $C$-short loop is trivial.
Redenoting, we may assume that every $C$-short loop in $S$ is trivial. This property holds for each covering of $S$.

To construct a covering $\widetilde{S}\rightarrow S$ satisfying 2), we will work with each boundary component of $S$ separately and then take a pullback.
Let $R$ be a boundary component of $S$. Choose a point $x\in R$.
Let $\gamma$ be a curve in $S$ with endpoints in $R$ such that $l(\gamma)\leqslant C$. By the above assumption, we
may assume that $\gamma$ is simple. 

Suppose that $\gamma$ is not homotopic relative to the endpoints to a curve in $R$. 
By applying a homotopy if needed, we may assume that $\gamma$ intersects $R$ only in the endpoints of $\gamma$.
Let $\gamma'$ be the simple closed curve starting at $x$ such that $\gamma$ is a part of $\gamma'$ and
$\gamma'\setminus \gamma\subseteq R$.
Clearly, $l(\gamma')\leqslant C+C_0$,
where $C_0$ is the maximum of lengths of boundary components of $S$.

Let $\widehat{S}$ be the surface obtained from $S$ by gluing a cup along $R$.
Suppose that $\gamma'$ is trivial in $\widehat{S}$.
Since $\gamma'$ is simple, it bounds a disc in $\widehat{S}$. Considering two cases,
where this disc contains the cup or not, we deduce that $\gamma$ can be homotoped, relative to the endpoints, into $R$ in $S$, a contradiction.
Thus, $\gamma'$ is nontrivial in $\widehat{S}$.

Then there is a normal covering $\overline{\widehat{S}}\rightarrow \widehat{S}$, where all lifts of $\gamma'$ have the property that their endpoints lie on different lifts of the cup. Removing all lifts of the cap from $\overline{\widehat{S}}$,
we get a covering $\overline{S}\rightarrow S$, where each lift of $\gamma'$ has endpoints on different lifts of $R$.
Hence, each lift of $\gamma$ has endpoints on different lifts of $R$.

We construct such coverings for each homotopy class $[\gamma']\in \pi_1(S,x)$ such that the length of the geodesic representative of $[\gamma']$ does not exceed $C+C_0$ and the image of $[\gamma']$ in $\pi_1(\widehat{S},x)$ is nontrivial. In this way we construct coverings for each boundary component of $S$. The pullback of all such coverings together with the identity covering $S\rightarrow S$ satisfies 2). \hfill$\Box$





\section{Appendix B: Constructing coverings\\ with prescribed branching data}

Recall that if $S$ is a compact surface of genus $g$ with $n$ boundary components, then $\chi(S)=2-2g-n$ if $S$ is orientable and $\chi(S)=1-g-n$ if not. The following proposition was proved first by D.~Husemoller in~\cite{Hus} und reproved in~\cite{Edm}. It is not valid for $S$ with genus 0.
The corresponding problem is called {\it Hurwitz realizability problem for branched coverings of surfaces}, see~\cite{Hur}.  It is not solved yet in general, see~\cite{Edm,P1,P2,Pak}.

\begin{prop}\label{222}
Let $S$ be a compact orientable surface of genus $g(S)\geqslant 1$ with boundary components $R_1,\dots, R_n$.
Let $\widetilde{S}$ be a compact orientable surface with boundary components $R_{i,j}$ and associated natural numbers $d_{i,j}$, $i=1,\dots,n$, $j=1,\dots,k_i$.

Then there exists a covering $\widetilde{S}\rightarrow S$ of degree $d$ such that $R_{i,j}$ covers $R_i$ with degree $d_{i,j}$ if and only if the following conditions are satisfied.

\medskip

{\rm 1)} $\chi(\widetilde{S})=d\cdot \chi(S)$,

{\rm 2)} $d=\overset{k_i}{\underset{j=1}{\sum}}d_{i,j}$ for every $i=1,\dots ,n$.

\end{prop}

The following proposition was proved by A.L.~Edmonds, R.S.~Kulkarni and R.E.~Stong in~\cite[Proposition~5.2]{Edm}.

\begin{prop}\label{2a} Let $S$ be a compact orientable surface of genus $g(S)= 0$ with boundary components $R_1,\dots, R_n$, $n\geqslant 3$.
Let $\widetilde{S}$ be a compact orientable surface with boundary components $R_{i,j}$ and associated natural numbers $d_{i,j}$, $i=1,\dots,n$, $j=1,\dots,k_i$.
Suppose that $k_n=1$ and $d_{n,1}=d$.

Then there exists a covering $\widetilde{S}\rightarrow S$ of degree $d$ such that $R_{i,j}$ covers $R_i$ with degree $d_{i,j}$ if and only if the following conditions are satisfied.

\medskip

{\rm 1)} $\chi(\widetilde{S})=d\cdot \chi(S)$,

{\rm 2)} $d=\overset{k_i}{\underset{j=1}{\sum}}d_{i,j}$ for every $i=1,\dots ,n$.
\end{prop}


\medskip

{\it Notation.} We denote the branching data of the covering $\widetilde{S}\rightarrow S$ by
$$\bigl( R_1(d_{1,1},\dots ,d_{1,k_1}),\dots ,R_n(d_{n,1},\dots ,d_{n,k_n}) \bigr).$$

\medskip

\begin{lem}\label{cor1} Let $S$ be a compact orientable surface with boundary components $R_1,\dots, R_n$, and suppose that $n\geqslant 3$ if $g(S)= 0$.
Then there exists a covering $\widetilde{S}\rightarrow S$ of degree $d\in \{2,4\}$ such that $g(\widetilde{S})\geqslant 1$
and $\widetilde{S}$ has a boundary component $\widetilde{R}_1$ which covers $R_1$ with degree $1$.
\end{lem}


{\it Proof.} We consider three cases.

{\it Case 1.} Suppose that $g(S)\geqslant 1$.

Let $\widetilde{S}$ be a compact orientable surface of genus $g(\widetilde{S})=2\cdot g(S)-1$ and with $2n$ boundary components. Then $g(\widetilde{S})\geqslant 1$ and $\chi(\widetilde{S})=2\cdot \chi(S)$. By Proposition~\ref{222}, there exists a covering $\widetilde{S}\rightarrow S$ of degree $d=2$
with the branching data $$\bigl(R_1(1,1),\dots ,R_n(1,1)\bigr).$$

\medskip

{\it Case 2.} Suppose that $g(S)=0$ and $n$ is odd.

Let $\widetilde{S}$ be a compact orientable surface
of genus $g(\widetilde{S})=\frac{3n-7}{2}$ with $n+1$ boundary components.
Then $g(\widetilde{S})\geqslant 1$ and $\chi(\widetilde{S})=4\cdot \chi(S)$.
By Proposition~\ref{2a}, there exists a covering $\widetilde{S}\rightarrow S$ of degree $d=4$
with the branching data $$\bigl(R_1(1,3),\, R_2(4), \dots ,R_n(4)\bigr).$$

\medskip

{\it Case 3.} Suppose that $g(S)=0$ and $n$ is even.

Let $\widetilde{S}$ be a compact orientable surface
of genus $g(\widetilde{S})=\frac{3n-8}{2}$ with $n+2$ boundary components.
Then $g(\widetilde{S})\geqslant 1$ and $\chi(\widetilde{S})=4\cdot \chi(S)$.
By Proposition~\ref{2a}, there exists a covering $\widetilde{S}\rightarrow S$ of degree $d=4$
with the branching data $$\bigl(R_1(1,3),\, R_2(1,3), R_3(4), \dots ,R_n(4)\bigr).$$

\noindent
We see that in each case one boundary component of $\widetilde{S}$ covers $R_1$ with degree 1. \hfill $\Box$

\begin{lem}\label{corM} Let $S$ be a compact orientable surface different from a disc and with $n\geqslant 1$ boundary components.
Let $\phi:\widetilde{S}\rightarrow S$ be a covering of degree $d$. Then for every natural $M$ which is a multiple of $(4dn)!$, there exists
a covering $\psi:\widehat{S}\rightarrow S$ which has even degree and factors through $\phi$, and such that
the restriction of $\psi$ to any boundary component of $\widehat{S}$ has degree $M$.
\end{lem}

{\it Proof.} 
The lemma is trivial in case where $g(S)=0$ and $n=2$. So, we may assume that if $g(S)=0$, then $n\geqslant 3$.
We will find $\psi$ as the composition of three coverings $\widehat{S}\overset{\tau}{\rightarrow} S'\overset{\sigma}{\rightarrow} \widetilde{S}\overset{\phi}{\rightarrow} S$
for appropriate $\widehat{S}, S',\tau$ and $\sigma$.

By Lemma~\ref{cor1}, there exists a covering  $\sigma: S'\rightarrow \widetilde{S}$ of degree at most 4 and such that $g(S')\geqslant 1$.
Let $L_1,\dots ,L_m$ be all boundary components of $S'$ and suppose that the restrictions of $\phi\circ \sigma$ to these  boundary components  have degrees $d_1,\dots ,d_m$.
Then $d_1+\dots +d_m$ equals to the degree of $\phi\circ \sigma$ multiplied by $n$. Hence
$d_1+\dots +d_m\leqslant 4dn$.

Let $M$ be a multiple of $(4dn)!$. Then each $d_i$ is a divisor of $M$.

Let $\widehat{S}$ be a compact orientable surface with $2(d_1+\dots +d_m)$ boundary components and of genus
$$g(\widehat{S})=M\cdot (2g(S')+m-2)-(d_1+\dots+d_m)+1.$$
This number is nonnegative since $g(S')\geqslant 1$, $m\geqslant 1$, and because of the choice of~$M$.
For this surface we have $\chi(\widehat{S})=2M\cdot \chi(S')$.
By Proposition~\ref{222}, there exists a covering $\tau:\widehat{S}\rightarrow S'$ of degree $2M$ with the branching data
$$\bigl(L_1(\underbrace{M/d_1, M/d_1, \dots ,M/d_1}_{2d_1}),
\dots ,L_m(\underbrace{M/d_m, M/d_m, \dots ,M/d_m}_{2d_m})\bigr).$$


Then we can put $\psi:=\phi\circ \sigma\circ \tau$. The degree of $\psi$ is even since the degree of $\tau$ is even.
\hfill $\Box$


\begin{lem}\label{corM1} Let $S$ be a compact orientable surface with $n\geqslant 1$ boundary components $R_1,\dots,R_n$. Suppose that $n\geqslant 3$ if $g(S)=0$.
Let $\phi:\widetilde{S}\rightarrow S$ be a covering of even degree $d$ with $g(\widetilde{S})\geqslant 1$ and
with the branching data
$$\bigl( R_1(d_{1,1},\dots ,d_{1,k_1}),\dots ,R_n(d_{n,1},\dots ,d_{n,k_n}) \bigr),$$
where $d_{1,1}=1$.
Then for every natural $M$ which is a multiple of
$4\overset{n}{\underset{i=1}{\prod}}\overset{k_i}{\underset{j=1}{\prod}}d_{i,j}$, there exists
a covering $\psi:\widehat{S}\rightarrow S$  which factors through $\phi$ and
has the branching data


$$\bigl(R_1(\underbrace{1, 1, \dots ,1}_{M/2},\,M/2,\,\underbrace{M, M, \dots ,M}_{d-1}),\,\, R_2(\underbrace{M, M, \dots ,M}_{d}), \dots ,R_n(\underbrace{M, M, \dots ,M}_{d})\bigr).
$$
\end{lem}

\vspace*{2mm}

{\it Proof.} Let $L_{1,1},\dots, L_{1,k_1},\dots,L_{n,1},\dots, L_{n,k_n}$ be the boundary components of $\widetilde{S}$, where we assume that $L_{i,j}$ covers $R_i$ with degree $d_{i,j}$.

\vspace*{1mm}

Let $\widehat{S}$ be a compact orientable surface with $dn+M/2$ boundary components and of genus $$g(\widehat{S})=\frac{1}{2}\bigl(M[2d\cdot g(S)+d(n-2)-1/2]+2-dn\bigr).$$

The latter number is integer since $4|M$ and $d$ is even. One can prove that this number is nonnegative
by considering two cases: $g(S)=0$, $n\geqslant 3$ and $g(S)\geqslant 1$, $n\geqslant 1$, and by using the fact that $M\geqslant 4$. One can check that $\chi(\widehat{S})=M\cdot \chi(\widetilde{S})$ by using the fact that $\chi(\widetilde{S})=d\cdot \chi(S)$.
By Proposition~\ref{222}, there exists a covering $\theta:\widehat{S}\rightarrow \widetilde{S}$ of degree $M$ with the branching data $(\vec{L}_{1,1},\dots, \vec{L}_{1,k_1},\dots,\vec{L}_{n,1},\dots, \vec{L}_{n,k_n})$, where
$$\begin{array}{ll}
\vec{L}_{1,1}= &L_{1,1}(\underbrace{1, 1, \dots ,1}_{M/2},M/2),\vspace*{3mm}\\
\vec{L}_{i,j}= & L_{i,j}(\underbrace{M/d_{i,j}, M/d_{i,j}, \dots ,M/d_{i,j}}_{d_{i,j}})
\hspace*{4mm} {\text{\rm if}}\hspace*{2mm} (i,j)\neq (1,1).
\end{array}
$$

Then we can set $\psi:=\phi\circ\theta$. The covering $\psi$
has degree $dM$ and it has the desired branching data, since $d_{1,1}=1$ and $d=\overset{k_i}{\underset{j=1}{\sum}}d_{i,j}$ for every $i=1,\dots ,n$.
\hfill $\Box$

\begin{lem}\label{very technical}
Let $S$ be a compact orientable surface with $\chi(S)\leqslant -1$, with
a fixed hyperbolic metric, and with $n\geqslant 1$ boundary components $R_1,\dots ,R_n$ which are geodesics
with respect to this metric.
Let $C_0$ be the maximum of lengths of boundary components of $S$.
For any constant $C>C_0$, there exists an even number $M_0$ satisfying the following statement:

for each multiple $M$ of $M_0$, there exist coverings $\Theta:\widehat{S}\rightarrow S$ and $\theta_i:\widehat{S}_{i}\rightarrow S$, $i=1,\dots,n$, such that:
\medskip

\begin{enumerate}
\item[{\rm (a$_1$)}] The branching data of $\Theta$ is
$$\bigl(R_1(\underbrace{M, M, \dots ,M}_{c}),\, R_2(\underbrace{M, M, \dots ,M}_{c}), \dots ,R_n(\underbrace{M, M, \dots ,M}_{c})\bigr)
$$ for some $c$.

\item[{\rm (a$_2$)}] There are no nontrivial $C$-short loops in $\widehat{S}$.

\item[{\rm (a$_3$)}] Every $C$-short geodesic curve in $\widehat{S}$
with endpoints on a boundary component of $\widehat{S}$
lies in this component.

\medskip

\item[{\rm (b$_1$)}] There exists $d\in \mathbb{N}$ such that the branching data of each $\theta_i$
can be obtained from the $n$-tuple
$$\bigl(R_1(\underbrace{M, M, \dots ,M}_{d}),\, R_2(\underbrace{M, M, \dots ,M}_{d}), \dots ,R_n(\underbrace{M, M, \dots ,M}_{d})\bigr)
$$
by replacing its $i$-th term by
$R_i(\underbrace{1, 1, \dots ,1}_{M/2},\,M/2,\,\underbrace{M, M, \dots ,M}_{d-1})$.

\item[{\rm (b$_2$)}] Every $C$-short loop in $\widehat{S}_i$ is freely homotopic into a 1-boundary of $\widehat{S}_i$.

\item[{\rm (b$_3$)}] Every $C$-short geodesic curve in $\widehat{S}_i$ with endpoints on a 1-boundary of $\widehat{S}_i$
lies in this 1-boun\-dary.

\item[{\rm (b$_4$)}] The distance between any two 1-boundary components of $\widehat{S}_i$ is larger than~\!$C$.

\end{enumerate}

\medskip
\noindent

\end{lem}

{\it Proof.} Let $C_1$ be the minimum of lengths of boundary components of $S$.

a) By Lemma~\ref{7a}, there exists a covering $\phi:\widetilde{S}\rightarrow S$ such that $\widetilde{S}$ satisfies the conditions (a$_2$) and (a$_3$).
By Lemma~\ref{corM}, for each $M$ which is a multiple of a certain even number, there exists a covering $\Theta_{M}:\widehat{S}\rightarrow S$ which factors through $\phi$ and satisfies (a$_1$). Since $\Theta_M$ factors through $\phi$, the surface $\widehat{S}$ satisfies (a$_2$) and (a$_3$) as~well.


b)
For $i=1,\dots ,n$, let $\psi_i:\widetilde{S}_i\rightarrow S$ be the covering from Lemma~\ref{lemma7} constructed for
the boundary component $R_i$. In particular, $\widetilde{S}_i$ has a boundary component $\widetilde{R}_i$ which covers $R_i$ with degree 1.

By Lemma~\ref{cor1}, we can construct a finite covering $\delta_i:\widetilde{S}_i\,'\rightarrow \widetilde{S}_i$ of even degree and such that $g(\widetilde{S}_i\,')\geqslant 1$, and $\widetilde{S}_i\,'$ has a boundary component $\widetilde{R}_i\,'$ which covers $\widetilde{R}_i$ with degree~1.
Set $\phi_i:=\psi_i\circ \delta_i$. Then, the covering $\phi_i:\widetilde{S}_i\,'\rightarrow S_i$ has the following properties:

i) $\phi_i$ factors through $\psi_i$ and has even degree;

ii) $g(\widetilde{S}_i\,')\geqslant 1$;

iii) $\widetilde{S}_i\,'$ has a 1-boundary $\widetilde{R}_i\,'$.

\medskip
\noindent
By Lemma~\ref{corM1}, for every natural $M$ which is a multiple of a certain even number, there
exists a covering $\theta_i:\widehat{S}_i\rightarrow S$ which factors through $\phi_i$ and satisfies (b$_1$).
By the statement 4) of Lemma~\ref{lemma7}, $\theta_i$ also satisfies (b$_3$) and (b$_4$). Moreover, by this statement, each $C$-short geodesic loop in $\widehat{S}_i$ can be freely homotoped into a boundary component of $\widehat{S}_i$.
If we choose $M>2C/C_1$,
we can deduce from (b$_1$) that this component is a $1$-boundary. Thus, (b$_2$)
is satisfied.

Finally, we can take $M_0$ as the least common multiple of two values of $M$ considered in a) and b).\hfill $\Box$



\end{document}